\documentclass[12pt,a4paper]{article}
\usepackage{algorithm}
\usepackage{algpseudocode}
\usepackage{amsmath}
\usepackage{bm}
\usepackage{url}
\usepackage{array}
\usepackage{subcaption}
\usepackage{tabularray}
\UseTblrLibrary{booktabs}
\usepackage{booktabs}
\usepackage{rotating}
\usepackage{caption}
\usepackage{multirow}
\usepackage{amssymb}
\usepackage{soul}
\usepackage{hyperref}
\usepackage[normalem]{ulem}
\usepackage{graphicx}
\usepackage{natbib}
\usepackage[dvipsnames]{xcolor}

\title{Sparse Identification of Time Delay Systems via Pseudospectral Collocation\thanks{The authors are members of INdAM research group GNCS; DB is a member of UMI research group Modellistica socio-epidemiologica''. EB and DB were supported by the European Union - NextGenerationEU with the project MOdellistica Numerica e Data-driven per l'Innovazione sostenibile - MONDI'' (CUP: G25F21003390007). MT was supported by the Italian Ministry of University and Research (MUR) through a PhD grant PNRR DM351/22.}}

\author{Enrico Bozzo\textsuperscript{1,2}, Dimitri Breda\textsuperscript{1,3}, Muhammad Tanveer\textsuperscript{1,4}}

\date{}

\begin{document}
\maketitle

\begin{center}
\textsuperscript{1}CDLab – Computational Dynamics Laboratory, Department of Mathematics, Computer Science and Physics, University of Udine, Via delle Scienze 206, 33100 Udine, Italy\\
\textsuperscript{2}Email: enrico.bozzo@uniud.it\\
\textsuperscript{3}Email: dimitri.breda@uniud.it\\
\textsuperscript{4}Email: tanveer.muhammad@spes.uniud.it
\end{center}

\begin{abstract}                
We present a pragmatic approach to the sparse identification of nonlinear dynamics for systems with discrete delays. It relies on approximating the underlying delay model with a system of ordinary differential equations via pseudospectral collocation. To minimize the reconstruction error, the new strategy avoids optimizing all possible multiple unknown delays, identifying only the maximum one. The computational burden is thus greatly reduced, improving the performance of recent implementations that work directly on the delay system.
\end{abstract}
\noindent%
\vfill

\section{Introduction}
In the last decades, data-driven model discovery has emerged as a lively research field due to an increased availability of data and easy-to-access computational resources \citep{bk19}. Under this paradigm, techniques to recover from measurements of the governing equations of an underlying dynamical system have gained a prominent role. In particular, the use of Sparse Identification of Nonlinear Dynamics (SINDy) has spread rapidly from its introduction for Ordinary Differential Equations (ODEs) in  \cite{bpk16}. Nowadays extensions exist to treat more general classes of problems, from partial \citep{rud17} to stochastic \citep{bnc18} differential equations, to name just a couple. Let us soon recall that the backbone of SINDy consists in expressing the Right-Hand Side (RHS) of the underlying dynamics as a linear combination of functions from a chosen library, and sparsity follows since in general only a few such functions are necessary \citep{bpk16}.

Only recently extensions of SINDy to Delay Differential Equations (DDEs) have been investigated, yet limited to constant discrete delays. To the best of our knowledge, \cite{sandoz23,pec24,wu23,kop24} are the only available references on the subject (a preliminary version of \cite{pec24} was first presented at IFAC TDS 2022, \cite{pec22}). In particular, \cite{sandoz23,pec22} propose the natural idea of adapting SINDy to DDEs by including delayed evaluations of the library functions. Moreover, the case of an unknown delay is tackled by minimizing the reconstruction error of SINDy over a set of candidate delays, respectively by following a Brute Force (BF) approach in \cite{sandoz23} (i.e., evaluating all the candidates) and via Bayesian Optimization (BO) in \cite{pec24}. An extension to multiple delays and other unknown parameters has also been considered in \cite{pec24}, always via BO on sets of candidates. The problem of optimally selecting candidate delays has also been analyzed in \cite{kop24}, while \cite{wu23} proposes the use of parameterized library functions to reduce the dimensionality issue (of a library ideally containing all candidate terms) by optimizing the concerned parameters via Particle swarm-based Optimization (PO in the following).

In this work we follow on the same line of the above contributions, focusing on effectively adapting SINDy to the case of unknown delays, both in number and values. We resume from the basic approach of identifying the correct delays by {\it externally} minimizing the reconstruction error of SINDy, first improving by using PO instead of BF or BO in terms of the number of calls to SINDy by the external optimizer. Then, in order to overcome the problem of dealing with intermediate multiple delays (i.e., those besides the maximum one), we propose a novel approach based on reducing the underlying DDE to an ODE via pseudospectral collocation, following \cite{bdgsv16}. This leads to a pragmatic tool that asks to externally optimize only the maximum delay, thus resorting to univariate optimization rather than a demanding multivariate one necessary to optimize all the intermediate delays. As a partial drawback, this new approach does not recover an interpretable RHS in the case of multiple delays, but it rather returns a ``black-box'' ODE able to match the given samples and to simulate the underlying dynamics beyond the time span of the original data (anyway, note that in general a successful matching of the true trajectory by SINDy does not necessarily correspond to a correct reconstruction, and hence interpretation, of the RHS). As further contributions, we introduce a physics-informed error; allow for both uniform and random distribution of samples (with possible addition of noise); reconstruct delayed and collocated samples by linear interpolation of available data, so to avoid the constraint of uniformly spaced samples as for the existing approaches.

The work is structured as follows. In Section \ref{sec:ESINDy} we summarize the basic SINDy approach for DDEs, treating both known (Section \ref{sec:Etauknown}) and unknown (Section \ref{sec:Etauunknown}) delays. In Section \ref{sec:PSINDy} we present the novel SINDy approach, recalling the pseudospectral collocation (Section \ref{sec:PSC}) and presenting its application in the SINDy framework (Section \ref{sec:extlib}). In Section \ref{sec:exp} we illustrate via experimental comparison on several DDEs the performance of all the mentioned approaches, in terms of both SINDy itself and the external optimization (via BF, BO and PO) of unknown delays/parameters. We conclude in Section \ref{sec:conc} with a short overview on possible future directions.
\section{Basic SINDy and extensions}\label{sec:ESINDy}
SINDy is a tool based on linear regression for data-driven model discovery originally introduced for ODEs
in \cite{bpk16}. In order to grasp how it works, let us consider the problem of identifying the RHS $f:\mathbb{R}^{n}\to\mathbb{R}^{n}$ defining an $n$-dimensional system of ODEs $x'=f(x)$, starting from the measurements of $x(t)$ at $m$ time instants $t_1,\ldots,t_m$ organized in a matrix $X\in\mathbb{R}^{m\times n}$ such that $x_{i,j}=x_{j}(t_{i})$, $i=1,\ldots,m$, $j=1,\ldots,n$. SINDy assumes that $f$ is well approximated by the linear combination of the functions $\theta_\ell$, $\ell=1,\ldots,p$, of a predefined library (containing say $1,x,x^2,\ldots,\sin x,\ldots$). Componentwise, this amounts to $f_{j}(x(t_i))\approx\sum_{\ell=1}^p \eta_\ell(x(t_i))\xi_{\ell,j}$, $i=1,\ldots,m$, $j=1,\ldots,n$.
Since $f(x(t_i))=x'(t_i)$, once derivative samples have been collected or approximated from the available measurements to form a matrix $X'\in\mathbb{R}^{m\times n}$, the above equations can be gathered in the linear system
\begin{equation}\label{SINDyODE}
X'=\Theta(X)\Xi,
\end{equation}
where $\Theta(X)\in\mathbb{R}^{m\times p}$ represents the library functions evaluated at the samples and the unknown $\Xi\in\mathbb{R}^{p\times n}$ has components $\xi_{\ell,j}$, $\ell=1,\ldots,p$, $j=1\ldots,n$. As a rule $m>p$, so that \eqref{SINDyODE} is overdetermined and typically solved columnwise for $\Xi$ by minimizing $\|X'-\Theta(X)\Xi\|_2+\lambda\|\Xi\|_1$.
The regularization through $\lambda>0$ enhances sparsity in $\Xi$ since frequently $f$ is the combination of just a few functions \citep{bpk16}.
Next, summarizing from \cite{sandoz23,pec24}, we address the case of DDEs with discrete delays, assumed to be known in Section \ref{sec:Etauknown} and unknown in Section \ref{sec:Etauunknown}.

\subsection{Extension to DDEs with known delays: E-SINDy}\label{sec:Etauknown}
Consider a DDE with $k$ constant discrete delays
\begin{equation}\label{DDE}
x'(t)=f(x(t),x(t-\tau_1),\ldots,x(t-\tau_k))
\end{equation}
for some $f:\mathbb{R}^{n(k+1)}\to\mathbb{R}^{n}$. When both the number and the values of the delays are known, \cite{sandoz23,pec24} extend SINDy for identifying $f$ by integrating the library functions with delayed copies. Hence \eqref{SINDyODE} is replaced by $X'=\Theta(X,X_{\tau})\Xi(\tau)$,
where $X_{\tau}$ is a shortcut for $X_{\tau_{1}},\ldots,X_{\tau_{k}}$ with each $X_{\tau_{i}}$ containing the delayed samples $x(t-\tau_{i})$ and the notation $\Xi(\tau)$ highlights the dependence of the sparse solution $\Xi$ on the given vector of delay values $\tau:=(\tau_{1},\ldots,\tau_{k})^{T}$.
Note that in both works the delayed samples are assumed to be available among the measurements, owing to the hypothesis that each delay is a multiple of a uniform sampling time. Here we avoid this assumption by getting delayed samples via linear interpolation of the original measurements, thus allowing both irrational delays and randomly distributed samples.

Let us remark that the above extension of SINDy to DDEs in the class \eqref{DDE} assumes the available knowledge of both the number and the values of the concerned delays. Imagining that this uncommon case is driven by a certain amount of modeling expertise (w.r.t. the capability, e.g., of interpreting the available data), we name this approach {\it E-SINDy}, where E stands for ``Expert''. Note indeed that E-SINDy operates directly on the underlying DDE once assumed it has the form \eqref{DDE}. Moreover, it potentially aims at recovering the RHS $f$ in view of, e.g., model interpretation and simulation beyond the time span of available data. As for the latter, note again that a successful matching of the true trajectory cannot exclude a failure in recovering the RHS.

\subsection{Extension to DDEs with unknown delays}\label{sec:Etauunknown}
In the more realistic case of unknown delays, a possible approach first considers the map 
\begin{equation}\label{epstau}
\tau\mapsto\epsilon(\tau):=\|X'-\Theta(X,X_{\tau})\Xi(\tau)\|_{2}
\end{equation}
giving the reconstruction error of E-SINDy. Then, for the case of a single unknown delay, \cite{sandoz23} proposes BF to select the value of $\tau$ minimizing $\epsilon(\tau)$ over a finite discrete set of candidate values. In \cite{pec24} this approach is improved in terms of calls to E-SINDy by resorting to BO. Therein, a multivariate BO is also invoked to minimize \eqref{epstau} over multiple delays, as well as for other unknown parameters (due, e.g., to specific nonlinear functions in the library as for the Mackey-Glass equation, see Section \ref{sec:MG}). In general, we call ``external optimization'' the minimization of \eqref{epstau} -- opposite to the ``internal optimization'' that SINDy performs for sparse regression (here we adopt sequential thresholded least-squares as in \cite{bpk16}). As a final remark, in the following we replace \eqref{epstau} with
\begin{equation}\label{epstauPI}
\setlength\arraycolsep{0.1em}\begin{array}{rcl}
\epsilon(\tau;w)&:=&w_{1}\|X'-\Theta(X,X_{\tau})\Xi(\tau)\|_{2}\\[1mm]
&&+w_{2}\|X-\mathcal{S}_{{\rm DDE}}(\Theta(X,X_{\tau})\Xi(\tau))\|_{2},
\end{array}
\end{equation}
where $w=(w_{1},w_{2})^{T}$ is a couple of weights and $\mathcal{S}_{{\rm DDE}}(\Theta(X,X_{\tau})\Xi(\tau))$ indicates the trajectory of the DDE recovered by E-SINDy simulated via any available numerical routine (we use MATLAB's \texttt{dde23}).
\section{A pragmatic SINDy approach}\label{sec:PSINDy}
E-SINDy for unknown delays as presented in Section \ref{sec:Etauunknown} requires to externally optimize (via BF or BO) a prescribed number $\bar k$ of delays assuming $\bar k\geq k$ where $k$ is the true number of delays in \eqref{DDE}. If $\bar k>k$ and E-SINDy is successful, $\bar k-k$ coefficients in the returned sparse matrix $\Xi$ will be null. Multivariate optimization is computationally demanding, so we propose a new approach where only the maximum delay (say $\bar\tau:=\tau_{k}$) is externally optimized. As a partial drawback, it will be clear from the following section that one renounces to recover an interpretable RHS, yet obtaining a ``black-box'' ODE matching the original data and correctly simulating the expected trajectory beyond the measurements time span. Moreover, we adopt PO to further improve computational efficiency. We name this new approach {\it P-SINDy}, where P stands now for ``Pragmatic'', following the terms coined in \cite{bdls16}. Indeed, opposite to an ``expert'' approach that deals directly with a DDE exploiting some prior knowledge, this new methodology relies on first reducing the DDE to a finite-dimensional system of approximating ODEs via pseudospectral collocation as originally proposed in \cite{bdgsv16}.
The reduction procedure is consolidated, easy to implement and acting on ODEs goes back to the standard SINDy approach of \cite{bpk16}. Above all, P-SINDy does not require any prior knowledge on the number and values of possible multiple (intermediate) delays. Yet one could legitimately argue that the number of state variables is increased due to the discretization of the original infinite-dimensional state space of \eqref{DDE} (despite its finite ``physical'' dimension $n$). Actually, as we illustrate next, P-SINDy works on just $n$ variables as E-SINDy does.

\subsection{Pseudospectral collocation}\label{sec:PSC}
Consider a generic DDE Initial Value Problem (IVP)
\begin{equation}\label{IVPDDEF}
\left\{\setlength\arraycolsep{0.1em}\begin{array}{ll}
x'(t)=F(x_t),&\quad t \geq 0,\\[1mm]
x(\eta)=\varphi(\eta),&\quad\eta\in[-\bar\tau,0],
\end{array}\right.
\end{equation}
where $x_t(\eta):=x(t+\eta)$, $\eta\in[-\bar\tau,0]$, represents the state at time $t$ of the associated dynamical system on the state space $\mathcal{X}:=C([-\bar\tau,0],\mathbb{R}^{n})$, $\varphi\in\mathcal{X}$ and $F:\mathcal{X}\to\mathbb{R}^{n}$ is a smooth RHS (note that $F(\psi)=f(\psi(0),\psi(-\tau_{1}),\ldots,\psi(-\tau_{k}))$ gives \eqref{DDE}). \eqref{IVPDDEF} is equivalent to the abstract IVP in \(\mathcal{X}\)
\begin{equation*}
\left\{\setlength\arraycolsep{0.1em}\begin{array}{rcl}
u'(t)&=&\mathcal{A}(u(t)),\quad t\geq0,\\[1mm]
u(0)&=&\varphi,
\end{array}\right.
\end{equation*}
where $\mathcal{A}:\mathcal{D}(\mathcal{A})\subseteq\mathcal{X}\to\mathcal{X}$ given by $\mathcal{A}(\psi)=\psi'$ with domain $\mathcal{D}(\mathcal{A})=\{\psi\in\mathcal{X}:\psi'\in\mathcal{X}\text{ and }\psi'(0)=F(\psi)\}$
is the infinitesimal generator of the strongly continuous semigroup \(\{T(t)\}_{t \geq 0}\) of solution operators $T(t):\mathcal{X}\to\mathcal{X}$, $T(t)\varphi:=x_t$. The equivalence is given by $u(t)=x_t$ for $\varphi\in\mathcal{D}(\mathcal{A})$ and in a mild sense otherwise since $\mathcal{D}(\mathcal{A})$ is dense in $\mathcal{X}$ \citep{diekmann95}. To reduce the DDE in \eqref{IVPDDEF} to a system of ODEs, let $\eta_{i}:=\frac{\bar\tau}{2}\left(\cos\left(\frac{i\pi}{M}\right)-1 \right)$, $i=0,1,\ldots,M$,
be the $M+1$ Chebyshev extremal nodes in $[-\bar\tau,0]$. Correspondingly, let $\mathcal{X}_{M}:=\mathbb{R}^{n(M+1)}$ be the finite-dimensional counterpart of $\mathcal{X}$ and consider relevant restriction and prolongation operators
\begin{equation*}
\setlength\arraycolsep{0.1em}\begin{array}{ll}
R_{M}:\mathcal{X}\to\mathcal{X}_{M},&\quad R_{M}\psi:=(\psi(\eta_0),\psi(\eta_1),\ldots,\psi(\eta_M)),\\[1mm]
P_{M}:\mathcal{X}_{M}\to\mathcal{X},&\quad P_{M}\Psi:=\sum_{j=0}^{M}\ell_j\Psi_j
\end{array}
\end{equation*}
for $\{\ell_0,\ell_1,\ldots,\ell_M\}$ the Lagrange basis for $\{\eta_{0},\eta_{1},\ldots,\eta_{M}\}$. Then $\mathcal{A}$ is discretized by $\mathcal{A}_{M}:\mathcal{X}_{M}\to\mathcal{X}_{M}$ given as ${[\mathcal{A}_M(\Psi)]}_0=F(P_M\Psi)$ and ${[\mathcal{A}_M(\Psi)]}_{i}=[R_M(P_M\Psi)']_{i}$, $i=1,\ldots,M$.
Consequently, the IVP \eqref{IVPDDEF} is approximated by the IVP for the blockwise system of $M+1$ $n$-dimensional ODEs in $\mathcal{X}_{M}$
\begin{equation}\label{DDEFM}
\left\{\setlength\arraycolsep{0.1em}\begin{array}{rcll}
U_{0}'(t)&=&F(P_M U(t)),&\\[1mm]
U_{i}'(t)&=&D_M U(t),&\quad i=1,\ldots,M,\\[1mm]
U_{i}(0)&=&\varphi(\eta_{i}),&\quad i=0,1,\ldots,M,
\end{array}\right.
\end{equation}
where $U(t)=(U_{0}(t),U_{1}(t),\ldots,U_{M}(t))^{T}\in\mathcal{X}_{M}$ and $D_{M}\in\mathbb{R}^{nM\times n(M+1)}$ has $n\times n$ block entries $d_{i,j}:=\ell_{j}'(\eta_{i})I_{n}$, $i=1,\ldots,M$, $j=0,1,\ldots,M$ ($I_{n}$ is the identity on $\mathbb{R}^{n}$). Above, $U_{i}(t)$, $i=0,1,\ldots,M$, approximates $x_{t}(\eta_{i})=x(t+\eta_{i})$ for $t\geq0$. The interest is in $U_{0}(t)\approx x(t)$, which however depends on all $U$ via the interpolation polynomial $P_{M}U$.

\subsection{Extended collocation library: P-SINDy}\label{sec:extlib}
Most of the ODEs in \eqref{DDEFM} (those with block-index $i=1,\ldots,M$) comes from the differentiation action of $\mathcal{A}$, so they are basically independent of the specific DDE in \eqref{IVPDDEF} if not for the maximum delay $\bar\tau$ affecting the collocation nodes. The only ($n$-dimensional) ODE affected by the original RHS $F$ through the boundary condition in $\mathcal{D}(\mathcal{A})$ is the first one. As a consequence, we employ SINDy to recover only this ODE, via sparse regression on $X'=\Theta(X,X_{\eta})\Xi_{M}(\bar\tau)$,
where now $X_{\eta}$ is a shortcut for $X_{\eta_{1}},\ldots,X_{\eta_{M}}$ with each $X_{\eta_{i}}$ containing the collocated samples $x(t+\eta_{i})$ and the notation $\Xi_{M}(\bar\tau)$ highlights the dependence of $\Xi$ on $\bar\tau$
and on the collocation degree $M$. It is now clear that the resulting P-SINDy acts on the physical dimension $n$ as E-SINDy does, while it has no requirements concerning possible intermediate delays in $(0,\bar\tau)$.

Regarding reconstruction and validation, given that for P-SINDy the only unknown delay is the maximum one, we introduce
\begin{equation}\label{epstauPIM}
\setlength\arraycolsep{0.1em}\begin{array}{rcl}
\epsilon_{M}(\bar\tau;w)&:=&w_{1}\|X'-\Theta(X,X_{\eta})\Xi_{M}(\bar\tau)\|_{2}\\[1mm]
&&+w_{2}\|X-\mathcal{S}_{{\rm ODE}}(\Theta(X,X_{\eta})\Xi_{M}(\bar\tau))\|_{2}
\end{array}
\end{equation}
accordingly to \eqref{epstauPI}. Now $\mathcal{S}_{{\rm ODE}}(\Theta(X,X_{\eta})\Xi_{M}(\bar\tau))$ indicates the trajectory of the {\it full} system of ODEs in \eqref{DDEFM} recovered by P-SINDy as far as the first block is concerned and then completed with the known differentiation part, as a whole simulated via any available numerical routine for ODEs (we use MATLAB's \texttt{ode45}). The term ``pragmatic'' is thus further justified, given that tools for ODEs are way more available and consolidated than those for DDEs. Let us remark again that while E-SINDy can in principle return an interpretable RHS thanks to the direct exploitation of the form \eqref{DDE}, P-SINDy does not (unless $k=1$). Nevertheless, from a ``pragmatic'' point of view we claim that the latter does not represent a true drawback, given that in general SINDy can give a good reconstruction error (and hence satisfactory trajectory matching) even in the presence of a loose RHS reconstruction. Finally, we underline that the external optimization of unknown delays for E-SINDy works on $k$ variables, while for P-SINDy it works on just one variable, thus greatly enhancing the overall performance as experimentally confirmed next.
\section{Numerical experiments}\label{sec:exp}
We compare the performance of both E-SINDy and P-SINDy on four different DDE models, together with BF, BO and PO as external optimizers of unknown delays (or parameters).
All the tests were run on a Windows 11 OS (CPU 5GHz, RAM 16Gb) by using our MATLAB implementations for SINDy (codes available at \url{http://cdlab.uniud.it/software}, MATLAB version R2024a) and MATLAB built-in BO and PO optimizers, respectively \texttt{bayesopt.m}
from the statistics and machine learning toolbox and \texttt{particleswarm.m}
from the global optimization toolbox. We worked on $m$ uniformly spaced samples obtained by integrating \eqref{IVPDDEF} with initial function $\varphi$ with MATLAB's \texttt{dde23} 
on the time window $[0,T]$, using the first 60\% portion for training and the remaining for validation via simulation. In the results (figures and tables) the letter E refers to E-SINDy and the letter P refers to P-SINDy (used with $M=5$, $10$ and $15$, denoted respectively by P5, P10 and P15). Values for $m$, $\varphi$ and $T$ will be specified for each model in the relevant section, together with the reference values of the parameters to be recovered by SINDy, the reference values for the unknown delays/parameters to be optimized externally, as well as the features of the adopted SINDy library (polynomial degree and possible non-polynomial terms).

Several tests have been performed beyond those presented next to, investigate the role of the weights $w=(w_{1},w_{2})^{T}$ in \eqref{epstauPI} and \eqref{epstauPIM}, as well as that of the samples' distribution and of possible additional noise. We did not appreciate any particular effect of these features on the overall results, so that in what follows we used uniformly spaced samples, absence of noise and $w_{1}=w_{2}=1$. We remark that the uniformly spaced samples do not constrain the values of the unknown delays, as linear interpolation for reconstructing delayed samples is adopted as anticipated.

\subsection{The delay logistic equation}\label{sec:DLE}
\begin{figure}
\centering
\includegraphics[width=0.8\linewidth,height=50mm]{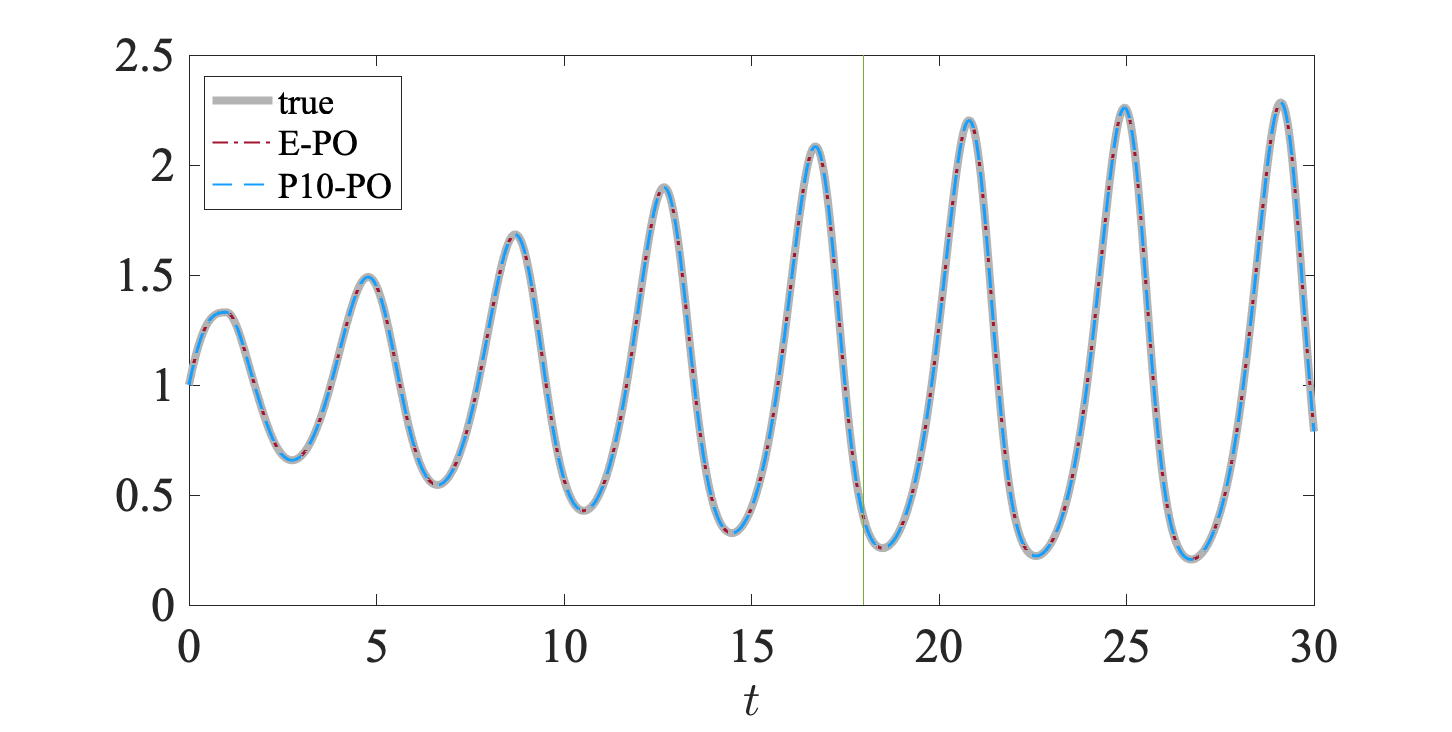}
\caption{Trajectory reconstruction for \eqref{DLE}, see text.}\label{fig:DLE1}
\end{figure}

\begin{figure}
\centering
\includegraphics[width=0.8\linewidth,height=50mm]{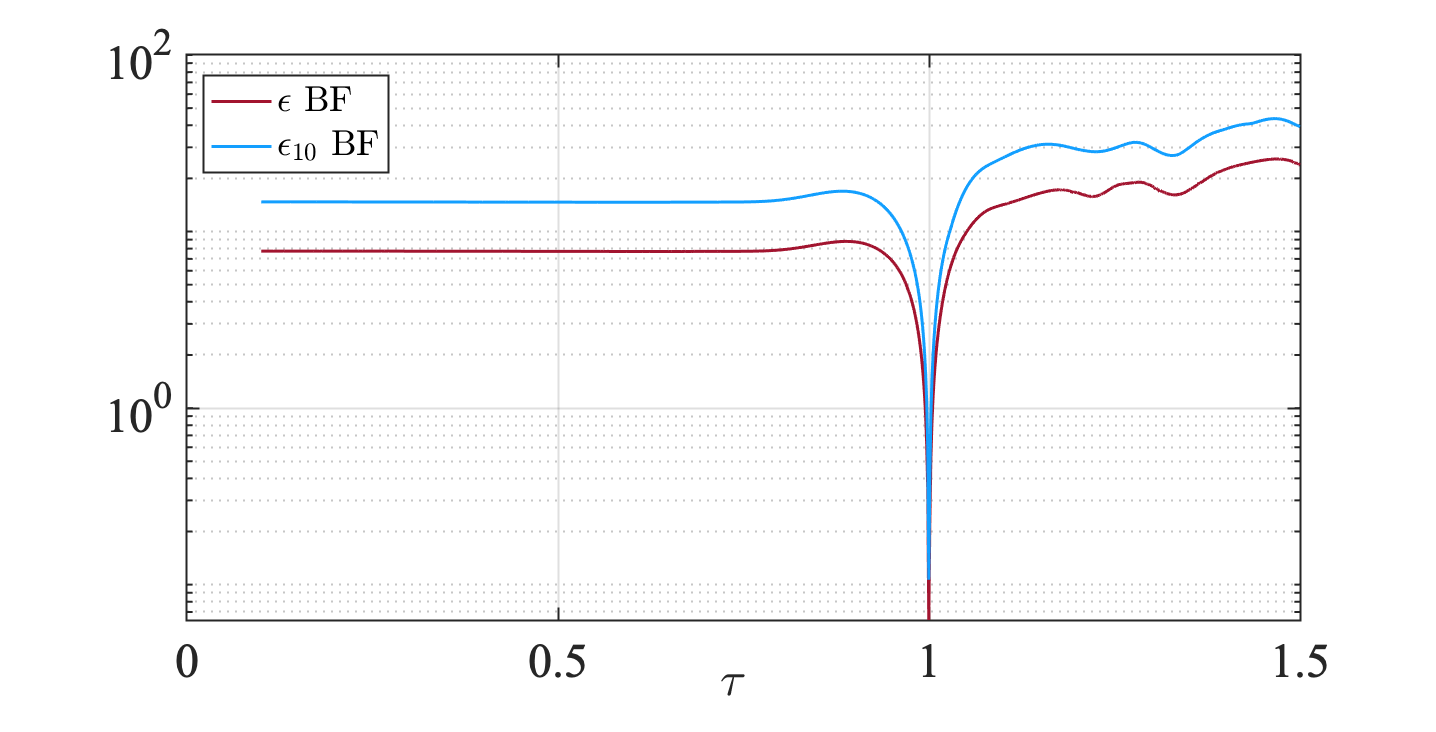}
\caption{Errors $\epsilon$ and $\epsilon_{10}$ for \eqref{DLE} with BF.}\label{fig:DLE2}
\end{figure}

\begin{figure}
\centering
\includegraphics[width=0.8\linewidth,height=50mm]{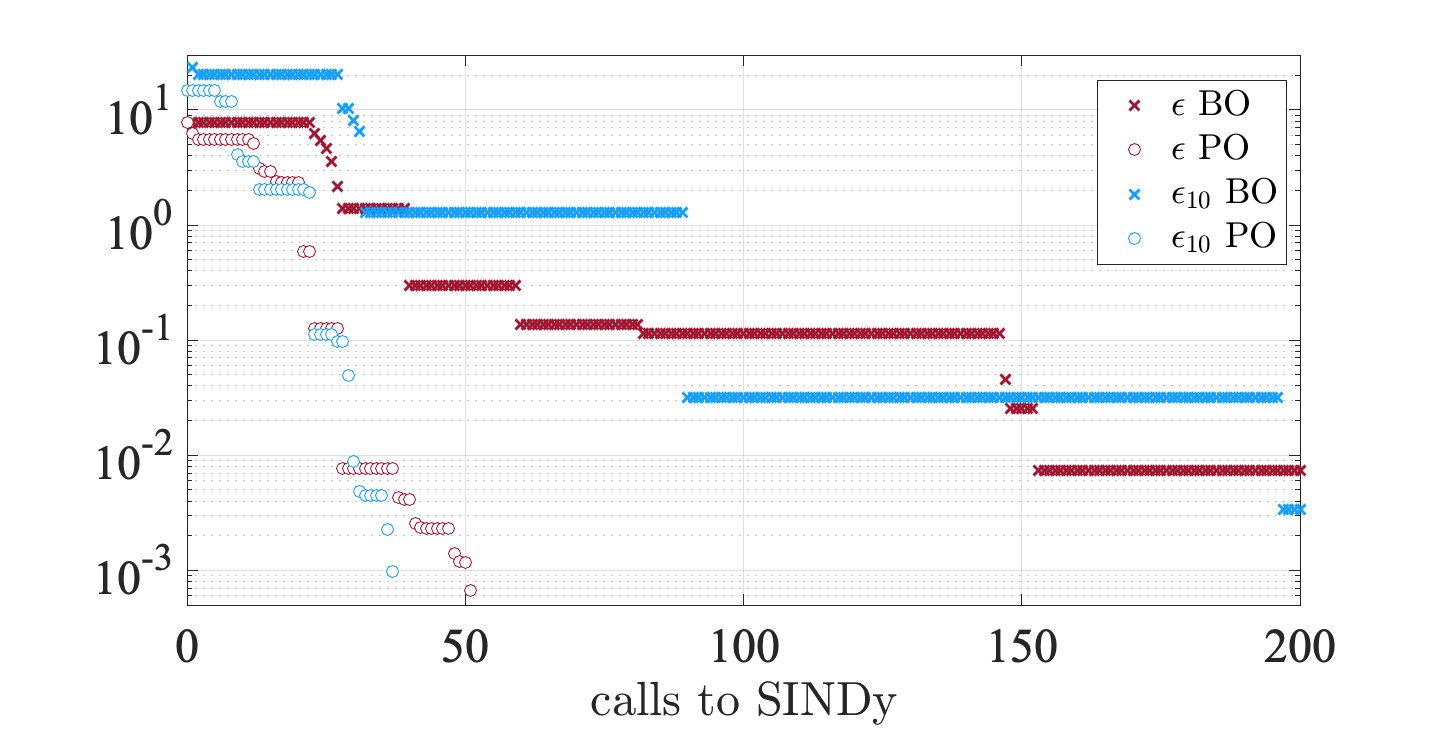}
\caption{Errors $\epsilon$ and $\epsilon_{10}$ for \eqref{DLE} with BO and PO.}\label{fig:DLE3}
\end{figure}

We consider the delay logistic equation \citep{hutch48}
\begin{equation}\label{DLE}
x'(t)=rx(t)(1-x(t-\tau))
\end{equation}
with reference parameter $r=1.8$ and a single unknown delay of true value $\tau=1$, using $m=100$, $T=30$ and $\varphi(\eta)=\cos(\eta)$. Figure \ref{fig:DLE1} shows relevant trajectories; Figure \ref{fig:DLE2} shows the reconstruction errors \eqref{epstauPI} and \eqref{epstauPIM} obtained via BF on $1\,000$ uniform candidate delays in $[0.1,1.5]$; Figure \ref{fig:DLE3} shows the behavior of the same errors when minimized via BO and PO vs the number of calls to SINDy. All the external optimizers returned the correct value of the true delay, yet with slightly different accuracy, Table \ref{tab:DLE}. Note that the accuracy of BF is dictated by the cardinality of the candidate set; BO stops at a fixed number of calls to SINDy given in advance;
PO stops when a desired accuracy is reached. Consequently, we fixed this accuracy to $10^{-3}$, chose the candidate set for BF accordingly and fixed the number of evaluations in BO in order to safely reach the minimum possible error. Both E-SINDy and P-SINDy gave back the correct sparse vector $\Xi$ with accuracy similar to that of the optimized value of $\tau$, Table \ref{tab:DLE}, using a polynomial library of degree 2. Correspondingly, the relevant trajectories are indistinguishable, both for training and validation, separated by the vertical thin green line in Figure \ref{fig:DLE1}. Therein, besides E-PO, we show only P10-PO as there is no remarkable difference in using $M=5$, $10$ or $15$, which means that a low collocation degree is sufficient to get a reasonable accuracy in such a data-driven context. Relevant numbers of calls to SINDy and CPU times are collected in Table \ref{tab:comparison1}, Section \ref{sec:comparison}, and commented therein together with the similar outcome for the other DDEs.
\begin{table}[htp]
\caption{Values and error returned for \eqref{DLE}.}\label{tab:DLE}
\centering
\begin{tabular}{llrrr}
\toprule
SINDy&optimizer&$r$&$\tau$ &$\epsilon$ or $\epsilon_{M}$\\
\midrule
E&BF&1.8000& 0.9997&6.0243e-5
\\
P5&BF&1.8001&0.9997&2.4373e-5\\
P10&BF&1.8000&0.9997&1.2565e-6\\
P15&BF&1.8000&0.9997&8.8321e-6\\
\midrule
E&BO&1.8001&0.9998&5.6233e-3\\
P5&BO&1.7999&1.0001&4.4588e-3\\
P10&BO&1.8000&1.0000&3.4849e-6\\
P15&BO&1.8000&1.0000&6.3592e-7\\
\midrule
E&PO&1.7999&1.0000&8.1706e-5\\
P5&PO&1.8100&1.0000&2.9017e-5\\
P10&PO&1.8000&1.0000&5.1410e-7\\
P15&PO&1.8000&1.0000&4.0167e-8\\
\bottomrule
\end{tabular}
\end{table}

\subsection{The Mackey-Glass equation}\label{sec:MG}
\begin{figure}
\centering
\includegraphics[width=0.8\linewidth,height=50mm]{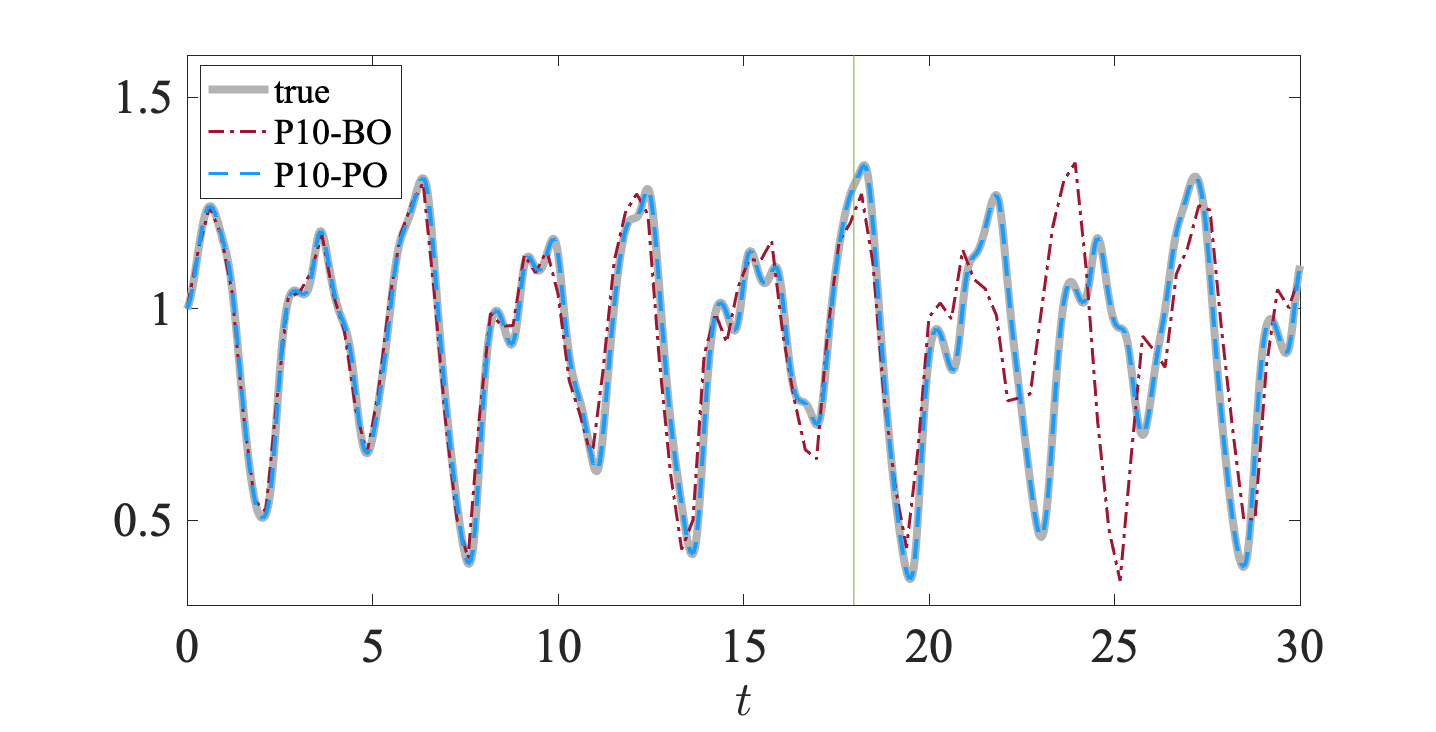}
\caption{Trajectory reconstruction for \eqref{MG}, see text.}\label{fig:MG}
\end{figure}

We now consider the Mackey-Glass equation \citep{mg77}
\begin{equation}\label{MG}
x'(t)=\beta\frac{x(t-\tau)}{1+x(t-\tau)^\alpha}-\gamma x(t)
\end{equation}
with reference parameters $\beta=4$ and $\gamma=2$ and unknown delay and exponent with reference values $\tau=1$ and $\alpha=9.6$, using again $m=100$, $T=30$ and $\varphi(\eta)=\cos(\eta)$. The polynomial library has degree 2 and includes the rational term $1/(1+X_{\tau}^{\alpha})$. In Figure \ref{fig:MG} we show the trajectories obtained with P10-BO and P10-PO, for which $\tau$ and $\alpha$ are simultaneously optimized. The relevant values are collected in Table \ref{tab:comparison2} of Section \ref{sec:comparison}: the final accuracy of BO is in general worse than that of PO, as it is evident from the resulting trajectories. The same can be said for E-SINDy instead of P-SINDy, but we avoid to plot the corresponding trajectories in favor of clarity. Comments on the number of calls to SINDy and CPU times are left to Section \ref{sec:comparison}, where for BF we used $100$ candidate delays in $[0.1,2]$ and $100$ candidate exponents in $[0.1,20]$.

\subsection{A scalar DDE with two delays}\label{sec:2tau}
\begin{figure}
\centering
\includegraphics[width=0.8\linewidth,height=50mm]{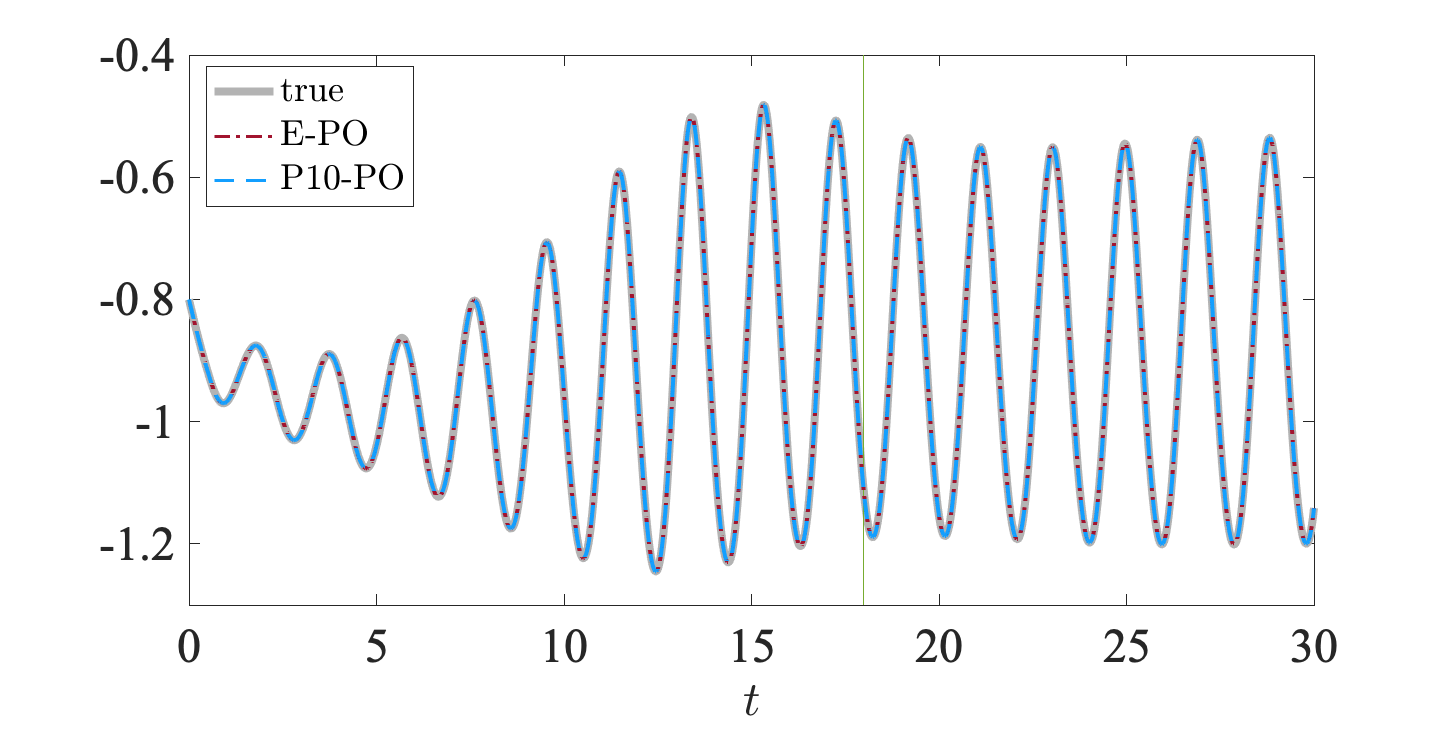}
\caption{Trajectory reconstruction for \eqref{2tau}, see text.}\label{fig:2tau}
\end{figure}

The DDE \citep{jmao21}
\begin{equation}\label{2tau}
x'(t)=a_2 x^2(t-\tau_1)+a_3 x^3(t-\tau_2)
\end{equation}
is tested with reference parameters $a_{2}=a_{3}=-1$ and unknown delays with reference values $\tau_{1}=0.65$ and $\tau_{2}=1.2$, using always $m=100$, $T=30$ and $\varphi(\eta)=\cos(\eta)$. The polynomial library has degree 3. In Figure \ref{fig:2tau} we show the trajectories reconstructed with E-PO and P1O-PO, which are again indistinguishable (as those obtained with BO, omitted for clarity). Comments on the number of calls to SINDy and CPU times are left to Section \ref{sec:comparison}, where for BF we used $100$ candidate delays for both $\tau_{1}$ and $\tau_{2}$, respectively in $[0.1,1]$  and in $[0.5,1.5]$. This choice amounts to $10\,000$ calls to SINDy for E-BF, as E-SINDy requires to optimize both $\tau_{1}$ and $\tau_{2}$. For P-SINDY we used just $1\,000$ candidate values for $\tau_{2}$, ensuring the same accuracy of PO.

\subsection{A delayed R\"ossler system}\label{sec:rossler}
\begin{figure}
\centering
\includegraphics[width=0.8\linewidth,height=50mm]{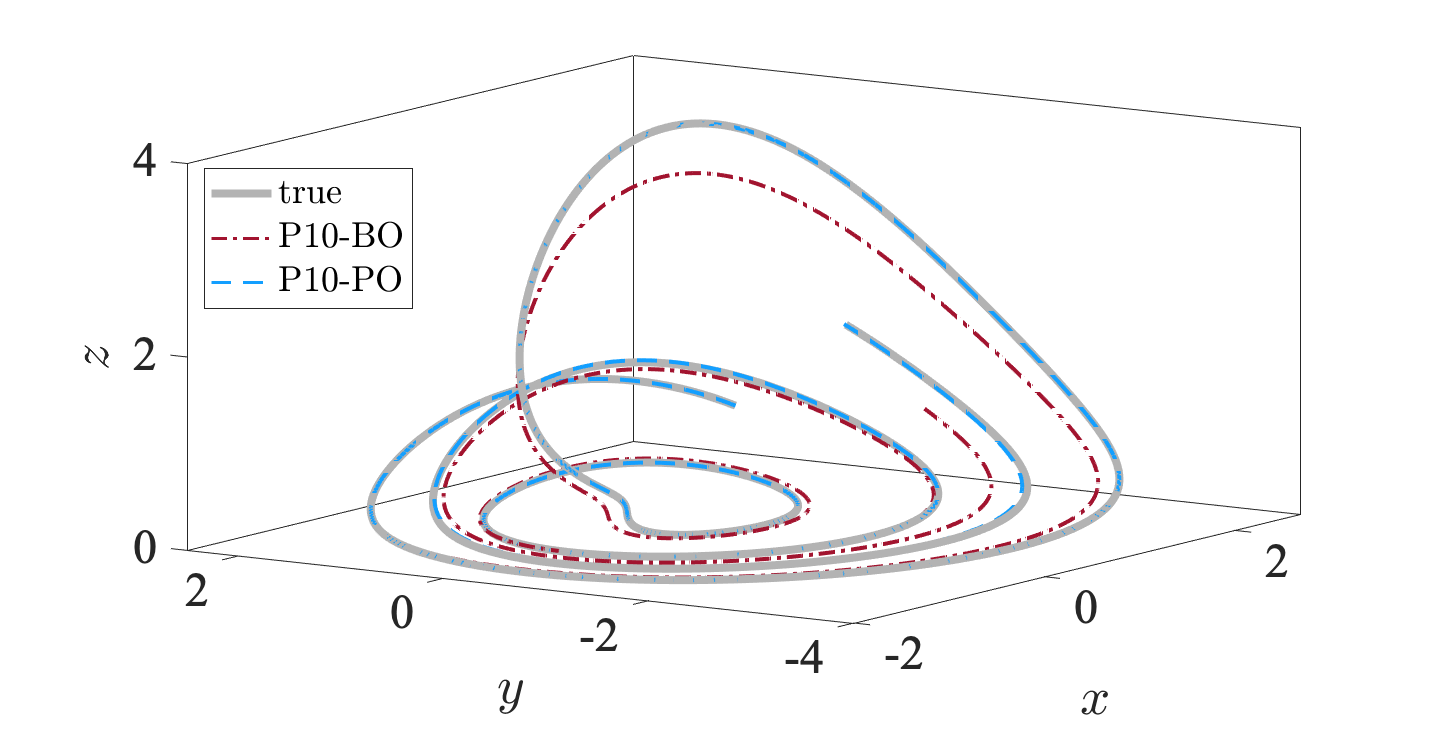}
\caption{Orbit reconstruction for \eqref{rossler}, see text.}\label{fig:rossler}
\end{figure}

The R\"ossler system with two delays \citep{wu23}
\begin{equation}\label{rossler}
\left\{\setlength\arraycolsep{0.1em}\begin{array}{rcl}
x'(t)&=&-y(t)-z(t)+\alpha_1 x(t-\tau_1)+\alpha_2 x(t-\tau_2),\\[1mm]
y'(t)&=&x(t)+\beta_1 y(t),\\[1mm]
z'(t)&=&\beta_2+z(t)(x(t)-1)
\end{array}
\right.
\end{equation}
is tested with reference parameters $\alpha_{1}=0.2$, $\alpha_{2}=0.5$, $\beta_{1}=\beta_{2}=0.2$ and unknown delays with reference values $\tau_{1}=1.5$ and $\tau_{2}=2$, using $m=100$, $T=30$ and $\varphi(\eta)=(1.5,0.4,0.9)^{T}$. The polynomial library has degree 2. The candidate delays are obtained by uniformly sampling $[0.1,2]$ for both $\tau_{1}$ and $\tau_{2}$ with 100 candidates. Relevant orbits in the physical state space $\mathbb{R}^{3}$ are depicted in Figure \ref{fig:rossler}, obtained with P10-BO and P10-PO. Note that, as for \eqref{MG}, PO performs better than BO, as the latter does not reach the same accuracy when optimizing the unknown delays. For further comments see Section \ref{sec:comparison}.

\subsection{General comparison and discussion}\label{sec:comparison}
In this section we summarize all the experiments, additionally collecting for all the SINDy methods (E, P5, P10, P15) and all the external optimizers (BF, BO, PO) the relevant number of calls to SINDy and CPU times in Table \ref{tab:comparison1}, as well as the reconstructed values of the unknown delays/parameters in Table \ref{tab:comparison2} (made exception for BF which gave the worst accuracy). The CPU times are given in seconds unless the total computation required more than an hour (in which case we simply write ``hrs''). The latter case occurred typically when using E-BF on two delays, as the cost scales w.r.t. the product of the cardinalities of the candidate sets and both delays are optimized (while P optimizes only the maximum one). The same happened also for P-SINDy for \eqref{MG} due to optimizing both $\tau$ and $\alpha$. In general one can observe that the cost of P-SINDy increases with the collocation degree $M$, as it is reasonable to expect. Good results are anyway obtained in general already with $M=5$, so that compared to E-PO, P5-PO is outperforming in terms of calls to SINDy, still balancing or even improving the CPU time. Only for \eqref{MG} we noted that $M=5$ was not enough to get a good trajectory reconstruction, most probably due to the possible presence of chaos (yet $M=10$ revealed sufficient, and still P10-PO required less evaluations in a comparable amount of time).

Summarizing, P-SINDy with PO improves the overall performance of both E-SINDy and BO. This superiority,  already evident for only one intermediate delay, further scales with the number of unknown delays.
\begin{table}[htbp]
\caption{Calls to SINDy and CPU times [s].}
\label{tab:comparison1}
\centering
\begin{tabular}{llrrrrrr}
\toprule
DDE & SINDy & \multicolumn{2}{c}{BF} & \multicolumn{2}{c}{BO} & \multicolumn{2}{c}{PO} \\
\cmidrule(r){3-4} \cmidrule(r){5-6} \cmidrule(r){7-8}
& & calls & CPU& calls & CPU & calls & CPU \\
\midrule
\eqref{DLE} & E & 1\,000 & 204 & 300 & 494& 308 & 67 \\
& P5 & 1\,000 & 217 & 300 & 438 & 241 & 70 \\
& P10 & 1\,000 & 449 & 300 & 485 & 179 & 100 \\
& P15 & 1\,000 & 1\,032 & 300 & 692 & 152 & 163 \\ \midrule
\eqref{MG}  & E & 10\,000 & hrs & 300 & 561 & 284 & 92 \\
& P5 & 10\,000 & hrs & 300 & 682 & 153 & 94 \\
& P10 & 10\,000 & hrs & 300 & 1\,080 & 51 & 102 \\
& P15 & 10\,000 & hrs & 300 & 1\,728 & 245 & 837 \\ \midrule
\eqref{2tau} & E & 10\,000 & hrs & 300 & 569 & 644 & 225 \\
& P5 & 1\,000 & 658 & 300 & 515 & 173 & 189 \\
& P10 & 1\,000 & 1\,406 & 300 & 1\,325 & 188 & 390 \\
& P15 & 1\,000 & 2\,697 & 300 & 2\,143 & 288 & 1\,111 \\ \midrule
\eqref{rossler} & E & 10\,000 & hrs & 300 & 576 & 220 & 197 \\
& P5 & 1\,000 & 1\,359 & 300 & 750 & 36 & 112 \\
& P10 & 1\,000 & 2\,130 & 300 & 1\,310 & 41 & 194 \\
& P15 & 1\,000 & hrs & 300 & 2\,376 & 31 & 214 \\
\bottomrule
\end{tabular}
\end{table}
\begin{table}[htbp]
\caption{Optimized unknown values.}
\label{tab:comparison2}
\centering
\begin{tabular}{llrrrrrr}
\toprule
DDE & SINDy & \multicolumn{2}{c}{BO} & \multicolumn{2}{c}{PO} \\
\cmidrule(r){3-4} \cmidrule(r){5-6} \cmidrule(r){7-8}
& & $\tau$ & $\alpha$&$\tau$&$\alpha$ \\
\midrule
\eqref{MG}  & E  & 0.9996 & 9.4969 & 1.0000 & 9.6001 \\
& P5  & 1.0502 & 9.2328 & 0.9994 & 9.4421 \\
& P10  & 1.0001 & 9.6134 & 1.0000 & 9.6001 \\
& P15& 0.9999 & 9.6001 & 1.0000 & 9.6000 \\ \midrule
& & $\tau_{1}$ & $\tau_{2}$ & $\tau_{1}$ & $\tau_{2}$ \\
\midrule
\eqref{2tau} & E  & 0.6478 & 1.2050 & 0.6500 & 1.2000 \\
& P5  & - & 1.2001 & - & 1.2000 \\
& P10  & - & 1.2000 & - & 1.2000 \\
& P15  & - & 1.2000 & - & 1.2000 \\ \midrule
\eqref{rossler} & E &1.5912 & 1.9684 & 1.5000 & 2.0000 \\
& P5  & - & 2.0000 & - & 2.0000 \\
& P10  & - & 2.0000 & - & 2.0000 \\
& P15  & - & 2.0000 & - & 2.0000 \\
\bottomrule
\end{tabular}
\end{table}
\section{Conclusions}\label{sec:conc}
We introduced the Pragmatic Sparse Identification of Nonlinear Dynamics (P-SINDy), a novel approach for the sparse identification of time-delay systems with discrete delays using pseudospectral collocation, which implicitly assumes to work with an ODE approximating the underlying DDE model. The results show that P-SINDy effectively handles unknown multiple delays by optimizing only the maximum one, thus significantly reducing computational demands. The method's ability to predict system dynamics beyond the available data using black-box simulations is particularly suitable when model knowledge is limited.

For future developments, we aim at further improving P-SINDy by integrating a piecewise pseudospectral collocation, with the consequent possibility of optimizing the collocation mesh. This is seen as a first step towards the challenge of identifying {\it distributed} delay terms via quadrature. Indeed, the latter would create several multiple intermediate delays (at the quadrature nodes), which can be efficiently handled by P-SINDy, while it would be prohibitive for standard SINDy implementations (as E-SINDy). These improvements would broaden the use of SINDy to more complex dynamical systems, including fundamental classes of models for structured population dynamics, driven in general by coupled delay and renewal equations \citep{dgg07}.

\end{document}